\documentclass[12pt]{article}

\usepackage{amsmath,amsthm,amssymb,amsfonts}
\usepackage{latexsym}
\def\pmb#1{\setbox0=\hbox{$#1$}%
\kern-.025em\copy0\kern-\wd0
\kern.05em\copy0\kern-\wd0
\kern-.025em\raise.0433em\box0}

\newfont{\lie}{eufm10 at 12pt}
\newfont{\field}{msbm10 at 11pt}
\newcommand{\R}{\mbox{\field \symbol{82}}}      
\newcommand{\h}{\mbox{\lie h}}                  
\newcommand{\Hyper}{\mbox{\field \symbol{72}}}  

\newcommand{\scr}[1]{\mbox{\scriptsize $#1$}}

\newcommand{\sm}[1]{\mbox{\small $#1$}}
\newcommand{\la}[1]{\mbox{\large $#1$}}
\newcommand{\La}[1]{\mbox{\Large $#1$}}

\newcommand{\Hu}[1]{\mbox{\Huge $#1$}}

\newcommand{\lnab}[1]{\la{\nabla}_{\!#1}}

\newcommand{\ra}{\rightarrow}
\newcommand{\non}{\nonumber}
\newcommand{\ha}{\scr{\frac{1}{2}}}
\newcommand{\al}{\alpha}

\newcommand{\Gf}{\Gamma_f}
\newtheorem{Lm}{Lemma}
\newtheorem{Pp}{Proposition}
\newtheorem{Th}{Theorem}

\begin{document}
\setcounter{page}{1}
\title{{Spacelike Graphs  with Parallel  Mean Curvature}}
\author{Isabel M.C.\ Salavessa}
\date{}
\protect\footnotetext{\!\!\!\!\!\!\!\!\!\!\!\!
 Partially supported by FCT through
POCI/MAT/60671/2004 and the Plurianual of CFIF.}
\maketitle
{\footnotesize \noindent
 Centro de F\'{\i}sica das Interac\c{c}\~{o}es Fundamentais,
Instituto Superior
T\'{e}cnico, Technical University of Lisbon,
Edif\'{\i}cio Ci\^{e}ncia,
Piso 3, Av.\ Rovisco Pais, 1049-001 Lisboa, Portugal;\\
e-mail: isabel.salavessa@mail.ist.utl.pt}\\[4mm]
{\small {\bf Abstract:} We consider spacelike graphs  $\Gf$
of simple products $(M\times N, g\times -h)$ 
where $(M,g)$ and $(N,h)$ are Riemannian
manifolds and  $f:M\ra N$ is a smooth map.
Under the condition of the Cheeger constant of $M$ to be zero and some
condition on the second fundamental form at infinity, 
 we conclude that
if $\Gf\subset M\times N$ has  parallel mean curvature $H$
then $H=0$.
This holds trivially if $M$ is closed. 
 If $M$ is the $m$-hyperbolic space then
for any constant $c$, we describe an explicit foliation
of $ \Hyper^m\times \R$ by hypersurfaces with constant mean curvature
$c$.
\markright{\sl\hfill Salavessa \hfill}
\section{Introduction}
\setcounter{Th}{0}
\setcounter{Pp}{0}
\setcounter{Cr} {0}
\setcounter{Lm} {0}
\setcounter{equation} {0}
The problem of estimating the mean curvature of a surface of
$\R^3$ described by  a graph of a function
$f:\R^2\ra \R$ was first introduced in 1955 by E.\ Heinz \cite{[H]}.
He proved that if $f$ is defined on the disc $x^2+ y^2<R^2$ and
the mean curvature satisfies $\|H\|\geq c>0$, where $c$ is a constant, 
then $R\leq \frac{1}{c}$. So, if $f$ is defined in all $\R^2$ and $\|H\|$
is constant, then $H=0$. Ten years later this problem was extended and 
solved   for the case of a map $f:\R^m \ra \R$ by Chern \cite{[Ch]} and
independently, by Flanders \cite{[F]}. In 1986, Jim Eells suggested
to the author a generalization of this problem in her
Ph.D thesis (\cite{[S1]}, \cite{[S2]}). We recall the formulation of
the problem.

Let $(M,g)$ and $(N,h)$ be two Riemannian manifolds of dimension $m$ and
$n$ respectively, and $f:M\ra N$ a smooth map. The graph of $f$,
$\Gf:=\{(p,f(p)):p\in M\}$ is a submanifold of $M\times N$ of dimension
$m$. We take on $M\times N$ the product metric $g\times h$,
 and on $\Gf$, the
induced one $\tilde{g}$. Let $H$ denote the mean curvature vector
of $\Gf$.
On $M$  it is defined the Cheeger constant
\[ \h (M)=\inf_D \frac{A(\partial D)}{V(D)}\]
where $D$ ranges over all open submanifolds of $M$ with compact
closure in $M$ and smooth boundary (see e.g. \cite{[Cha]}), and
$A(\partial D)$ and  $V(D)$ are respectively the area of $\partial D$ and the
volume of $D$, with respect to the metric $g$.
This constant is zero, if, for example, M is a closed manifold, or
if $(M,g)$ is a simple Riemannian manifold, i.e., there exists
a diffeomorphism $\phi:(M,g)\ra (\R^{m}, <,>)$ onto $\R^m$
such that $\lambda g \leq \phi^*<,>\leq \mu g$ for some positive constants
$\lambda, \mu$. Then we have got
\begin{Th} (\cite{[S1]},\cite{[S2]})
If $\Gf$ has parallel mean curvature with $c=\|H\|$, 
then for each oriented  compact domain 
$D\subset M$ we have the isoperimetric inequality
\[c\leq \frac{1}{m}\frac{A(\partial D)}{V(D)}.\]
Thus, $c\leq \frac{1}{m} \, \h (M)$. In particular if $(M,g)$ has
zero Cheeger constant then $\Gf$ is in fact a minimal
 submanifold of  $M\times N$.
\end{Th}
\noindent
In case $N$ is oriented one dimensional 
with unit vector field $"1"$, we do not
need parallel mean curvature to obtain a formula
\begin{equation}
m\, \langle H,\nu \rangle_{(g\times h)}=div_g\La{(}
\frac{\nabla f}{\sqrt{1+\|\nabla f\|^2}}\La{)}
\end{equation}
where $\nu=\frac{(-\nabla f, 1)}{\sqrt{1+\|\nabla f\|^2}}$ 
is  a unit normal to $\Gf$.
This led to a more general result:
\begin{Th} (\cite{[S1]},\cite{[S2]}) If $N$ is oriented of dimension one
 and $f:M\ra N$ is  any map, then
\begin{equation}
\min_{\bar{D}}\|H\|\leq  \frac{1}{m}\frac{A(\partial D)}{V(D)}.\\[1mm]
\end{equation}
\end{Th}
\noindent
This generalizes the inequality of Heinz-Chern-Flanders to
graphs of functions $f:M\ra \R$.
We note that it is not possible to relax the
assumption of $H$ to be constant to $0\leq H\leq C$, where $C$ is 
a constant, 
without further assumptions, to conclude minimality. In fact
we have the following example:
Set $f:\R^2\ra \R$ given by $f(x,y)= e^x$. 
Then
$ 0= \lim_{x\ra \pm \infty}H<  \ha div \la{(}
\frac{\nabla f}{\sqrt{1+\|\nabla f\|^2}}\la{)}
=\ha {e^x}{(1+ e^{2x})^{-\frac{3}{2}}} \leq C.$\\[3mm]
A more difficult kind of problem is the so-called Bernstein-type 
problems, that amounts to determine  geometric conditions
 to conclude that a minimal submanifold must be totally geodesic.
Recently Rosenberg \cite{[Ro]} obtained a Bernstein type result
for entire minimal graphs in $M^2\times \R$.
Al\'{\i}as, Dajczer and Ripoll have obtained in \cite{ [A-D-R]} 
a Bernstein-type result  for surfaces
in an ambient space a three-dimensional Riemannian manifold
endowed with a homothetic Killing field, that includes the case
of \cite{[Ro]}:
\begin{Th} (\cite{[A-D-R]}) Let $M^2$ be a complete surface with
Gauss curvature $K_M\geq 0$.\\[2mm]
$(i)$ Any entire constant mean curvature graph in $M^2\times \R$
 is totally geodesic.\\[1mm]
$(ii)$ If, in addition, $K_M(q)>0$ at some point $q\in M$,
then the graph is a slice.
\end{Th}
\noindent
The proof  is inspired in the ideas of Chern \cite{[Ch2]}
of a proof
of a Bernstein theorem in case $m=2$, and consists on
computing the Laplacian of  the support 
function $\Theta=\langle \nu, e\rangle $,
where $\nu$ is a globally defined unit normal vector field of
the normal bundle, and $e$ is a constant unit vector field tangent
to the factor $\R$.
The assumption
 of $K_M\geq 0$ is necessary, since in \cite{[D-H]}
and \cite{[N-R]} it is shown
the existence of non-trivial entire minimal graphs when $M^2=\Hyper^2$ is
the hyperbolic plane. \\[3mm]
\noindent
In the case $M=\Hyper^m$ the Cheeger constant is $(m-1)$. We have
an example constructed by the author  in \cite{[S1]} 
of a  graph in $\Hyper^m\times \R$ with non-zero constant mean curvature:
\begin{Pp}   Consider the hyperbolic space $\Hyper^m= (B^m,g)$ where 
$B^m$ is the unit open disk in $\R^m$ with centre $0$
and $g$ is the complete metric given by $g=4|dx|^2/(1-|x|^2)^2$, 
of constant sectional curvature equal to $-1$.
Let $c\in [1-m, m-1]$ and $f_c:\Hyper^m\ra \R$ defined by:
\[f_c(x)=\Hu{\int}_0^{r(x)}
\frac{\frac{c}{(\sinh r)^{m-1}}\int_0^r(\sinh t)^{m-1}dt}
{\sqrt{1-\left(\frac{c}{(\sinh r)^{m-1}}\int_0^r(\sinh t)^{m-1}dt\right)^2}}
dr, \]
where $r(x)=\log \left(\frac{1+|x|}{1-|x|}\right)$ 
is the distance function in $\Hyper^m$ to $0$. 
Then $f_c$ is smooth on all $\Hyper^m$, and for
each $d\in \R$, $\Gamma_{f_c +d}\subset
\Hyper^m\times \R$ has constant mean curvature given
 by $\|H\|=\frac{|c|}{m}$.
In the particular case $m=2$ and $c=1$, $f_c$ can be written as
\[f_c(x)=\int^{r(x)}_0\sqrt{\ha(\cosh r -1)}dr =\frac{2}{\sqrt{1-|x|^2}}-2\]
\end{Pp}
\noindent
In \cite{[S2]} we only 
give a brief explanation that this example exists in \cite {[S1]}. 
So we will give in the Appendix (section 3)
the proof of Proposition 1.1, that  reproduces the proof 
in \cite{[S1]}. Moreover, a slightly modified proof of this one gives 
a proof of Proposition 1.3 for the Lorentzian case. We also note the
following:
\begin{Pp} For each $c\in [1-m,m-1]$, $\{\Gamma_{f_c+d}(x): x\in \Hyper^m,
d\in \R\}$ defines a foliation of $\Hyper^m\times \R$ by hypersurfaces
with constant mean curvature $c$.
\end{Pp}
\noindent
\em Remark. \em If we fix $d$ and let $c$ to  vary, then
we also have a foliation  of $\Hyper^m\times \R$,
 on a neighbourhood
of $\Hyper^m\times \{d\}$,
by hypersurfaces
with constant mean curvature $c$, with $c$ varying on each leaf.
This also holds for the Lorentzian case of Proposition 1.3.
The author would like to thank the referee for pointing out this
interesting  detail.\\[3mm]

In this note we study the case of $M\times N$ is endowed with the
pseudo-Riemannian metric $g\times -h$. 
We abusively still call "minimal"  submanifolds, the ones that satisfy
 $H=0$. Note that a  graph $\Gf$ is  spacelike 
iff $f^*h\leq bg$ with $b:M\ra \R$ a continuous locally 
Lipschitz function
 satisfying $0\leq b(p)<1$, 
$\forall p\in M$. If $N$ is one-dimensional then
$b=\|\nabla f\|$.  In section 2 we  will prove the following:

\begin{Th}
Let $\Gf$ be a spacelike graph with parallel mean curvature with $\|H\|=
\sqrt{|\langle H,H\rangle|}=|c|$.
Then $\|\lnab{}df\|\geq \sqrt{m}|c|
(1-b)^2$, with equality iff $\lnab{}df=c=0$. 
Furthermore, if $\h (M)=0$ and if $\|\lnab{}df\|=O((1-b)^2)$,
then $\Gf$ is minimal. This is the case of $M$ compact.
\end{Th}
\noindent
Assume $N$ is oriented one-dimensional 
and $f:M\ra N$ defines a spacelike graph.
Then $\nu=-\frac{(\nabla f, 1)}{\sqrt{1-\|\nabla f\|^2}}$ is a unit 
timelike vector field that spans the normal bundle, and
defines a timelike direction. $H$ is
future directed  if $H=-\|H\|\nu$, with $\|H\|\geq 0$.
\indent
\begin{Th} Assume $N$ is oriented one-dimensional
 and $f:M\ra N$ defines a spacelike graph with future directed mean 
curvature.
On a compact domain $D$, let 
$b_D=\max_{\bar{D}}\|\nabla f\|=\max_{\bar{D}}b$.
Then
\begin{equation}
\min_{\bar{D}}\|H\|\leq \frac{1}{m}\frac{b_D}{\sqrt{1-b^2_D}}
\frac{A(\partial D)}{V(D)}.
\end{equation}
In particular, if $(1)$ or $(2)$ below holds:\\[1mm]
$(1)$  $\Gf$ has constant mean curvature, $\h (M)=0$, 
and $b\leq C<1$  for some
constant $C$;\\[1mm]
$(2)$ $|H|$ and ${b}/
({\sqrt{1-b^2}})$ are both integrable on $M$;\\[2mm]
 then $\Gf$ is a minimal spacelike
hypersurface. This is the case of $M$ compact.
\end{Th}
\noindent
If $b$ is not bounded by a constant $C<1$
or $\h(M)\neq 0$, we have an example,
very similar to the one of Proposition 1.1, except on a  sign 
in some term of the denominator.
\begin{Pp} 
Let $c$ be any constant and $f_c:\Hyper^m\ra \R$ defined by:
\[f_c(x)=\Hu{\int}_0^{r(x)}
\frac{\frac{c}{(\sinh r)^{m-1}}\int_0^r(\sinh t)^{m-1}dt}
{\sqrt{1+\left(\frac{c}{(\sinh r)^{m-1}}
\int_0^r(\sinh t)^{m-1}dt\right)^2}} dr, \]
where $r(x)=\log \left(\frac{1+|x|}{1-|x|}\right)$
is the distance function
in $\Hyper^m$ to $0$. Then $f_c$ is smooth on all 
$\Hyper^m$, and for each $d\in \R$, $\Gamma_{f_c+d}\subset
\Hyper^m\times \R$ is a spacelike graph 
with constant mean curvature given by 
$\langle H, \nu \rangle =\frac{c}{m}$. Furthermore,
$\{\Gamma_{f_c+d}(x):~ x\in \Hyper^m, d\in \R\}$ defines a foliation
of $\Hyper^m\times \R$ by hypersurfaces with constant mean curvature $c$.
\end{Pp}
\noindent
Examples of spacelike constant mean curvature  $H=c$ hypersurfaces
of $\R^{n+1}_1$ are the hyperboloids, i.e. the graph of
$f(x)=\sqrt{\frac{k^2}{m^2}\frac{1}{c^2}+\sum_{i=1}^kx_i^2 }$,
for $k=1,\ldots, n$. If $k=n$ this example and the ones of
Propositions 1.1 and 1.3 are described as constant mean curvature
 graphs of a function
$f:(M,g)\ra \R$ of the form 
$f(x)=\phi(r(x))$,  where $r(x)$ is the distance function in $(M,g)$
to a fixed point and $\phi:\R\ra \R$ is
a smooth function. Such $f$ are in fact smooth maps because $r^2$ is so,
and  $\phi$ (unique for a chosen constant $c$) can be expressed in 
terms of $r^2$.
\\[3mm]
It has been a relevant problem in General Relativity the study of the
 existence and uniqueness of space-like hypersurfaces 
with constant mean curvature in globally hyperbolic
connected Lorentzian manifolds having a compact Cauchy surface (GHLCS),
and the existence of foliations by such hypersurfaces. Here
we are  treating only the case $M\times N$ with a simple product
$g\times -h$. The metric of a GHLCS
is conformally equivalent to a
a warped product metric. 
For example, 
if $(M,g)$ is  closed  and $\R$ is endowed with the metric
$dt^2$ and  $\al:M\ra \R$ is any positive smooth function
( the lapse function),
then spacelike graphs $\Gf$ in $(M\times \R, g -\al^2 dt^2)$ 
exist with prescribed
 mean curvature $H:M\ra \R$ for any function $H$ satisfying
$\int_M H\al Vol_M=0$. These graphs are unique up to a constant (i.e, if
$\Gf$ is a solution then $\Gamma_{f+d}$ is also a solution). This was 
proved by Akutagawa (\cite{[Aku]}) using the invertibility of the Laplace
operator for closed $M$. In particular, if $H$ is constant, then
the submanifold must be minimal. On the other hand
on Robertson-Walker
spacetimes the slice hypersurfaces have constant mean curvature, and
recently, Al\'{\i}ais and Montiel \cite{[A-M]} proved that under certain 
conditions on the warping function, these are the only closed examples.
Gerhardt \cite{[Ger]} proved that  GHLCS spaces
can be foliated by constant mean curvature hypersurfaces if the big bang
and the big crunch hypothesis is satisfied and if a time-like
convergence condition holds. 

\section{Spacelike graphs}
\setcounter{Th}{0}
\setcounter{Pp}{0}
\setcounter{Cr} {0}
\setcounter{Lm} {0}
\setcounter{equation} {0}
Now we take in the product $M\times N$ the pseudo-Riemannian metric
$g\times -h$. If $f:M\ra N$ then we denote by
\[\begin{array}{cccc}
\Gamma_f:& M&\ra &(M\times N, g\times h)\\
 & p &\ra &(p,f(p))
\end{array}
\]
and identify the set $\Gamma_f$ with the embedding $\Gf$, and
let $\tilde{g}=\Gf^*(g\times -h)=g-f^*h$.
Assume that $f$ satisfies $h(df(X),df(X))<g(X,X)$. Then $\tilde{g}$
is a Riemannian metric of $\Gf$, that is $\Gf$ is a space-like
submanifold of $(M\times N, g\times -h)$. Let $H$ denote the mean
curvature of $\Gf$. Note that $H$ is a time-like vector.
\\[2mm]
Let $X_i$ a local o.n. frame of $(M,g)$ and $\tilde{g}_{ij}
=g(X_i,X_j)-h(df(X_i),df(X_j))$.   Set
\begin{eqnarray}
W=trace_{g-f^*h}(\lnab{}df)\in C^{\infty}(f^{-1}TN),\\[-1mm]
Z=\sum_{ij}\tilde{g}^{ij}h(W,df(X_i))X_j\in C^{\infty}(TM)
\end{eqnarray}
The following formulas hold:
\begin{Lm} If $\Gf$ has parallel mean curvature, then:\\
(1) $mH=(Z,W+df(Z))=(0,W)^{\bot}$.\\
(2) $m^2c^2= div_g(Z)$, where $c^2=-\langle H,H\rangle_{(g\times -h)}$.
\end{Lm}
\noindent
\em Proof. \em The proof is very similar to the one of
lemmas 1,2 and 3 of
 \cite{[S2]}
with some adjustments on the sign of $h$. So we omit it.  \qed\\[4mm]
Taking $X_i$ a o.n. basis of $T_pM$ that diagonalizes $f^*h$, i.e,
$df(X_i)=\lambda_i e_i$, for $i\leq k$
where $e_i$ is an o.n. system of $T_{f(p)}N$,
 and $df(X_i)=0$ for $i\geq k+1$, we conclude that
$a g\leq f^*h\leq bg$, and
$\frac{1}{1-a}g\leq\tilde{g}^{-1}\leq\frac{1}{1-b}g$ where
$a=\inf_i \lambda_i^2$ is the smallest eigenvalue of $f^*h$ and $b
=\sup_i \lambda_i^2$ the largest.
If $N$ is one-dimensional  and $m\geq 2$, then $a=0$ and 
$b=\|\nabla f\|.$ If we reorder the
eigenvalues $b=\lambda^2_1\geq \lambda^2_2\geq \ldots\geq  
\lambda^2_n=a$, including repeated eigenvalues according 
their multiplicity, by the Weyl's perturbation theorem each 
$\lambda_i^2$ is a continuous
locally Lipschitz function. In particular
 $b:M\ra [0,1)$ is a continuous locally Lipschitz function.\\[2mm]
From (2.1)-(2.2) we conclude:
\begin{Lm}
$~\|Z\|\leq \frac{\sqrt{b}}{1-b}\|W\|,$ and
~$\|W\|\leq \frac{\sqrt{m}}{1-b}\|\lnab{} df\|$.
\end{Lm}
\noindent
\em Proof of Theorem 1.4. \em\\[1mm]
Let $\Gf$ be a spacelike graph with parallel mean curvature. 
Using Lemma 2.1
\begin{eqnarray*}
-m^2c^2&=&\langle (Z,W+df(Z)), (Z,W+df(Z))\rangle_{g\times -h}\\
&=& \|Z\|^2-\|W\|^2-2h(W, df(Z))-\|df(Z)\|^2~~~~~~
\end{eqnarray*}
Thus,  
\begin{eqnarray}
\|Z\|^2 &\leq& -m^2c^2 +\|W\|^2+2\|W\|\|df(Z)\|+ \|df(Z)\|^2\non\\
&\leq& -m^2c^2 +(\|W\|+ \sqrt{b}\|Z\|)^2\non\\
&\leq& -m^2c^2+ \frac{m}{(1-b)^4}\|\lnab{}df\|^2
\end{eqnarray}
what implies the first assertion.
If $\|\lnab{}df\|^2=m c^2(1-b)^4$, then from (2.3) $Z=0$. Consequently,
by lemma 2.1(b), $c=0$, and so $\lnab{}df=0$.
 Now denote by $\bar{n}$ the unit
outward of $\partial D$. By lemma 2.1(2), Stokes and  Lemma 2.2
\begin{eqnarray}
m^2c^2V(D)=
\int_{\partial D}g(Z, \bar{n})\leq \int_{\partial D}\|Z\|
&\leq& A(\partial D)\sup_{\bar{D}}\frac{\sqrt{mb}}{(1-b)^2}\|\lnab{}df\|
\end{eqnarray}
If $\|\lnab{}df\|=O((1-b)^2)$,
  there exist a constant $C>0$ s.t.\ $\|\lnab{}df\|\leq C(1-b)^2$.
Then, from (2.4), $m^2c^2 \leq C'\frac{A(\partial D)}{ V(D)}$ 
for some constant $C'$
and Theorem 1.4 is proved.\qed\\[5mm]
Now
let us now assume  $N$ is oriented of dimension one with global
 vector field $"1"$. 
If $p_i=h(df(X_i),1)$, then $\tilde{g}_{ij}=\delta_{ij}-p_ip_j$,
$\tilde{g}^{ij}=\delta_{ij}+\frac{p_ip_j}{(1-\|\nabla f\|^2)}$. 
Similarly to the Riemannian case \cite{[S1]}, we can obtain a formula:
\begin{Lm}
\begin{equation}
m\langle H, \nu\rangle=div_g\left(\frac{\nabla f}
{\sqrt{1-\|\nabla f\|^2}}\right)
\end{equation}
\end{Lm}
\noindent
\em Proof of Theorem 1.5. \em \\[2mm]
We obtain (1.3) by integration over $D$ 
of (2.5) and use Stokes. (1) is an immediate consequence of the
definition of $\h (M)$, and (2) is a consequence of the extended theorem
of Stokes due to Gaffney \cite{[Gaf]} applied to (2.5).\qed 
\section{Appendix}
\setcounter{Th}{0}
\setcounter{Pp}{0}
\setcounter{Cr} {0}
\setcounter{Lm} {0}
\setcounter{equation} {0}
\subsection{Proof of Propositions 1.1 and 1.2} 
First we note that 
if $f$ satisfies $c=(1.1)$ then it does so $f+d$, where $d$ is 
a constant.
The function $r(x)=\log \left(\frac{1+|x|}{1-|x|}\right) 
=2\tanh^{-1}(|x|)$ has the following properties: $\forall x\neq 0$,
$\nabla r= \frac{1-|x|^2}{2}\frac{x}{|x|}$, where the gradient
of $r$ is w.r.t. the metric $g$. Hence , $\|\nabla r\|_g=1$
and $\Delta r= (m-1)\coth r$. We observe that $r^2$ is smooth.
Let us write $f= \phi(r)$ with $\phi:\R^+_0\ra \R$. Then
$\nabla f= \phi'(r) \nabla r$, and  so (1.1) applied to $f$ becomes
equivalent to 
\begin{eqnarray*}
c&=& div\La{(}\frac{\nabla f}{\sqrt{1+\|\nabla f\|^2_g}}\La{)}
=div\La{(}\frac{\phi'(r)\nabla r}{\sqrt{1+(\phi'(r))^2}}\La{)}\\
&=& \frac{\phi'(r)\Delta r}{\sqrt{1+(\phi'(r))^2}}
-\frac{(\phi'(r))^2\phi''(r) \|\nabla r\|^2}
{(1 + (h'(r))^2)^{\frac{3}{2}}} + \frac{
\phi''(r) \|\nabla r\|^2_g}{\sqrt{1+(\phi'(r))^2}}
\end{eqnarray*}
Using the above properties of $r$ we get
\begin{eqnarray*}
\lefteqn{c(1+(\phi'(r))^2)^{\frac{3}{2}}=}\\
&=& (m-1)\coth r(\phi'(r))(1+(\phi'(r))^2)
-(\phi'(r))^2\phi''(r) + \phi''(r)(1+((\phi'(r))^2)\\
&=& (m-1)\coth r(\phi'(r))(1+(\phi'(r))^2)+ \phi''(r)
\end{eqnarray*}
With the substitution $w(r)=\phi'(r)$ the equation becomes
\begin{equation}
w'= c(1+w^2)^{\frac{3}{2}}-(m-1)\coth r\, w \,(1+w^2), ~~~~\forall r>0
\end{equation}
The next step is to reduce this differential equation to a linear one
through several changes of variables. First write (3.1) as
\[
\frac{w\, w'}{(1+w^2)^{\frac{3}{2}}} =cw 
-(m-1)\coth r\, \frac{w^2}{(1+w^2)^{\frac{1}{2}}}
\]
Let $y=\frac{1}{(1+w^2)^{\frac{1}{2}}}\in (0,1]$. Then
$w= \pm \frac{\sqrt{1-y^2}}{y}$. Assume first the sign $+$. Then
\[(3.1) ~~\Longleftrightarrow~~ 
-y'=c\frac{\sqrt{1-y^2}}{y}-(m-1)\coth r\, \frac{1-y^2}{y^2}\,y.\]
Thus,
$ -y y'= c\sqrt{1-y^2}-(m-1)\coth r\, (1-y^2).$
Let $v=y^2\in (0,1]$. Then
\[(3.1) ~~\Longleftrightarrow~~ -\ha \frac{v'}{\sqrt{1-v}}=c-(m-1)\coth r
\,\sqrt{1-v}.\]
Finally, let $u=\sqrt{1-v}\in [0,1)$. Hence
\begin{equation}
(3.1) ~~\Longleftrightarrow~~ u'=c - (m-1)\coth r \, u
\end{equation}
which equation is linear. Let us first suppose $c=1$. Then, the general
solution of (3.2) is given by
\begin{eqnarray*}
u(r)&=& e^{-\int_{r_0}^{r}(m-1)\coth t dt}\left( \int_{r_0}^r
e^{(m-1)\int_{r_0}^{s}(m-1)\coth t dt}ds + u_0\right)\\
&=& e^{-(m-1)(\log \sinh r-\log \sinh r_0)}\left(
\int_{r_0}^r e^{(m-1)(\log \sinh s-\log \sinh r_0)}ds + u_0\right)\\
&=& \frac{(\sinh r_0)^{m-1}}{(\sinh r)^{m-1}}\left(\frac{1}{
(\sinh r_0)^{m-1}}\int_{r_0}^r
(\sinh s)^{m-1}ds + u_0\right)\\
&=& \frac{1}{(\sinh r)^{m-1}}\int_{r_0}^r
(\sinh s)^{m-1}ds+ u_0\frac{(\sinh r_0)^{m-1}}{(\sinh r)^{m-1}}.
\end{eqnarray*}
Let us now put $r_0=u_0=0$. Then we have 
\begin{equation}
u(r)=\frac{1}{(\sinh r)^{m-1}}\int_{0}^r(\sinh s)^{m-1}ds, ~~~~~~\forall r>0
\end{equation}
Next we prove that $u\in [0,1)$ with $u(0)=0$, and, moreover, that
$\sup_{r\in (0,+\infty)}u(r)=\lim_{r\ra +\infty}u(r)=\frac{1}{m-1}$.
Obviously $u$ is positive and with l'Hospital rule,
\[u(0)=\lim_{r\ra 0}u(r)=\lim_{r\ra 0}\frac{(\sinh r)^{m-1}}
{(m-1)(\sinh r)^{m-2}\cosh r}=\lim_{r\ra 0}\frac{\tanh r}{(m-1)}=0.\]
If $u(r)$ attains a local maximum at some $r_0\in (0, +\infty)$, then
$u'(r_0)=0$. From (3.2) we have $u(r_0)=\frac{\tanh r_0}{m-1}$.
Thus, $u(r_0)<\frac{1}{m-1}\leq 1$. On the other hand, if there are no
local maxima, then, necessarily, $\sup_{r\in (0,+\infty)}u(r)=
\lim_{r\ra +\infty }u(r)$. So only we have to calculate this limit.
With partial integration
\begin{eqnarray*}
\lefteqn{\int_0^r(\sinh s)^{m-1}ds=}\\
&=& \la{[}\cosh s(\sinh s)^{m-2}\la{]}^r_0-(m-2)\int_0^r
\cosh^2 s(\sinh s)^{m-3}ds\\
&=& \cosh r(\sinh r)^{m-2}-(m-2)\int_0^r(1+\sinh^2s)(\sinh s)^{m-3}ds\\
&=& \cosh r(\sinh r)^{m-2}-(m-2)\int_0^r(\sinh s)^{m-3}ds -
(m-2)\int_0^r(\sinh s)^{m-1}ds.
\end{eqnarray*}
Thus $\int_0^r(\sinh s)^{m-1}ds=\frac{1}{m-1}\cosh r(\sinh r)^{m-2}
-\frac{m-2}{m-1}\int_0^r(\sinh s)^{m-2}ds$, and
\begin{eqnarray*}
\frac{\int_0^r(\sinh s)^{m-1}ds}{(\sinh s)^{m-1}}
&=& \frac{1}{m-1}\coth r - \frac{(m-2)}{(m-1)\sinh ^2r}
\frac{\int_0^r(\sinh s)^{m-3}ds}{(\sinh r)^{m-3}}.
\end{eqnarray*}
Since $\forall p$, $ \frac{\int_0^r(\sinh s)^{p}ds}{(\sinh s)^{p}}$
is a bounded function on $r\in [0, +\infty)$, we have
\[\lim_{r\ra +\infty}\frac{\int_0^r(\sinh s)^{m-1}ds}{(\sinh s)^{m-1}}
=\frac{1}{m-1}\lim_{r\ra +\infty}\coth r =\frac{1}{m-1}.\]
Therefore,
\begin{equation}
\sup_{r\in[0,+\infty)}u(r)=\frac{1}{m-1}
\end{equation}
which is not a maximum. So, $0\leq u(r)<\frac{1}{m-1}$, $\forall r\in
[0, +\infty)$ and $u(r)$ satisfies (3.2) for $c=1$. Let now $c$ be
any arbitrary constant. Then, the function $\tilde{u}(r)=c u(r)$ is a 
solution of (3.2), but we have to impose $\tilde{u}(r)\in [0,1)$. 
From (3.4) we conclude that $c$ must satisfy $0\leq c\leq m-1$. That is,
$\forall 0\leq c\leq m-1$, the function
\[ \tilde{u}(r)= c\frac{\int_0^r(\sinh s)^{m-1}ds}{(\sinh s)^{m-1}}\]
fulfills the condition specified in (3.2).
In terms of the original function $f$, we have $f$ given by the 
expression in the Prop.1.2. 
If we had chosen the sign $-$ for the expression of $w$ we would 
get in (3.2) a $-c$ instead $c$ and we would obtain $\tilde{u}$ 
with a change of sign,
or equivalently, the same expression as in the Proposition, with
$c\in[-m+1,0]$.
Obviously, $f$ is smooth on $\Hyper^m\sim \{0\}$. Let us now investigate
the behaviour of $f$ close to the origin. Near $t=0$ we have 
the following  Taylor expansions:
\begin{eqnarray*}
&&\sinh t = t + \sm{\frac{t^3}{6}}+ O(t^5)=
t(1+\sm{\frac{t^2}{6}} +O(t^4))\\[-1mm]
&&(1+t)^m= 1+mt +\theta(t^2)\\[-1mm]
&&\frac{1}{\sqrt{1+t}}=1-\sm{\frac{t}{2}}+\theta(t^2),~~~~
\frac{1}{1-t}=1+t+\theta(t^2)
\end{eqnarray*}
where $\theta(t)$ and $O(t^k)$ are analytic functions of the form
\[\theta(t^k)=\sum_{n\geq 0}\sm{\frac{a^{k+n}}{(k+n)!}}t^{k+n}~~~
O(t^k)=\sum_{n\geq 0}\sm{\frac{a^{k+2n}}{(k+2n)!}}t^{k+2n}\]  
Then we have $\frac{1}{\sqrt{1+t^2}}=1-\sm{\frac{t^2}{2}}+\theta(t^4)$,
$\frac{1}{1-t^2}=1+t^2+\theta(t^4)$, and 
\[
(\sinh t)^{m-1}= t^{m-1}(1+ \sm{\frac{t^2}{6}}+O(t^4))^{m-1}
=t^{m-1}(1+ \sm{\frac{(m-1)}{6}}t^2)+ O(t^{m+3}).
\]
Hence
\begin{eqnarray*}
\lefteqn{
\frac{1}{(\sinh s)^{m-1}}\int_0^s(\sinh t)^{m-1}dt =
\frac{\frac{s^m}{m}+\frac{m-1}{m-2}\frac{s^{m+2}}{6} + O(s^{m+4})}
{s^{m-1}(1+\frac{(m-1)}{6}s^2 +O(s^4))}=}\\
&=&\frac{\frac{s}{m}+\frac{m-1}{m-2}\frac{s^{3}}{6} + O(s^{5})}
{(1+\frac{(m-1)}{6}s^2 +O(s^4))}\\
&=& \left(\sm{\frac{s}{m}}
+\sm{\frac{m-1}{m-2}\frac{s^{3}}{6}} + O(s^{5})\right)
\left(1- s^2(\sm{\frac{m-1}{6}}+ O(s^2))+ O(s^4)\right)\\
&=& s\left(1-\sm{\frac{(m-1)}{6}}s^2\right)
\left(\sm{\frac{1}{m}}+ \sm{\frac{(m-1)}{(m+2)}
\frac{s^2}{6}}\right)+ O(s^5)\\
&=& \sm{\frac{s}{m}}\left(1- \sm{\frac{(m-1)}{(m+2)}\frac{s^2}{3}}
\right) + O(s^5)
\end{eqnarray*}
For $A$ close to zero,  $\frac{A}{\sqrt{1-A^2}}=A(1+\ha A^2)+ O(A^5)$.
Putting
\[ A=\frac{c}{(\sinh s)^{m-1}}\int_0^s(\sinh t)^{m-1}dt
=c\sm{\frac{s}{m}}\left(1- \sm{\frac{(m-1)}{(m+2)}\frac{s^2}{3}}
\right) + O(s^5)\]
we have $O(A^5)=O(s^5)$ and
\begin{eqnarray*}
\lefteqn{\frac{A}{\sqrt{1-A^2}}=}\\
&=&\La{(}s\sm{\frac{c}{m}}
\left(1- \sm{\frac{(m-1)}{(m+2)}\frac{s^2}{3}}\right) + O(s^5)\La{)}
\La{(}1+ \ha\la{(} \sm{\frac{cs}{m}}(1-\sm{\frac{(m-1)}{(m+2)}\frac{s^2}{3}})
+O(s^5)\la{)}^2\La{)} + O(s^5)\\
&=&s \sm{\frac{c}{m}}\left( 1+ s^2(\sm{\frac{c^2}{2m^2}}
-\sm{\frac{(m-1)}{3(m+2)}})\right)
+ O(s^5).
\end{eqnarray*}
Therefore
\begin{eqnarray*}
\int_0^r\frac{A}{\sqrt{1-A^2}} ds&=&{\frac{c}{m}\frac{r^2}{2}}
+{\frac{c}{m}\frac{r^4}{4}}\La{(} {\frac{c^2}{2m^2}}
-{\frac{(m-1)}{3(m+2)}}\La{)}
+O(r^6)={\frac{c}{m}\frac{r^2}{2}}+O(r^4)
\end{eqnarray*}
Consequently
\[f(x)=\int_0^{r(x)}\frac{A}{\sqrt{1-A^2}}\, ds = 
\frac{c}{m}\frac{r^2(x)}{2}+O(r^4(x)).\]
Since $r^2(x)$ is smooth on all $\Hyper^m$, we conclude that $f(x)$ 
is, too.\\
Finally, for each $x\neq 0$ fixed, the function $c\ra f_c(x)$ 
has non-zero derivative. Moreover,
$\Gamma_{f_c+d}(x)=\Gamma_{f_c+d'}(x')$ implies $x=x'$ and $d=d'$.
So we have two possible foliations, either varying $c$ or $d$.
Note that $O(r^4)$ also depends on $c$.
\qed.
\subsection{Proof of Proposition 1.3}
We solve $c=(2.5)$ for $f=\phi(r)$.
In this case we follow the previous proof, with the following
replacements:
\begin{eqnarray*}
c &=& div\left( \frac{\phi'( r) \,\nabla r}{\sqrt{1-(\phi'(r))^2}}
\right)\\
 w=\phi'(r),~~~|w|<1&& w'=c(1-w^2)^{\frac{3}{2}}-(m-1)\coth r \, 
w (1-w^2)\\
y=\frac{1}{\sqrt{(1-w^2)}}\in [1, +\infty) && v=y^2\in [1, +\infty)\\
u=\sqrt{v-1}\in [0,+\infty) && u'= c-(m-1)\coth r\, u
\end{eqnarray*}
Thus $u(r)$ is the same function as in the Riemannian case, but
now we do not have any restriction on the range of values of $u(r)$.
This implies we may choose first $u(r)$ as defined  in (3.3),
that corresponds to take $c=1$, and next  take $\tilde{u}=cu$
for any constant $c$ with no restrictions on the chosen $c$.
Finally, the proof that $f$ is smooth close the origin we use
$\frac{A}{\sqrt{1+A^2}}=A(1-\ha A^2)+ O(A^5)$ obtaining as well
\[f(x)=\int_0^{r(x)}\frac{A}{\sqrt{1+A^2}}ds =
\frac{c}{m}\frac{r^2(x)}{2}+O(r^4(x)).\]
and proving its smoothness.\\

We also note that the hyperboloid with $k=n$ is obtained in the same 
way, by taking $r(x)=\|x\|$ the Euclidean norm.
\qed


\begin{thebibliography}{29}
\baselineskip .3cm
\bibitem{[Aku]}\small{K.\ Akutagawa, \em A note on spacelike hypersurfaces 
with prescribed mean curvature in a spatially closed globally static 
Lorentzian manifold, \em  Mem.\ Fac.\ Sci.\ Kyushu Univ.\ Ser.\ A 
{\bf 40} (1986),  no. 2, 119--123.}
\bibitem{[A-D-R]}\small{L.J.\ Al\'{\i}as; M.\ Dajczer; J.\ Ripoll, \em A
Bernstein-type theorem for Riemannian manifolds with a Killing
field, \em preprint.}
\bibitem{[A-M]}\small{L.J.\ Al\'{\i}as; S.\ Montiel, 
\em Uniqueness of spacelike hypersurfaces with constant mean curvature in 
generalized Robertson-Walker spacetimes.\em  Differential geometry, Valencia, 
2001,  59--69, World Sci. Publishing, River Edge, NJ, 2002.}
\bibitem{[Cha]}\small{I.\ Chavel, \em  
Eigenvalues in Riemannian geometry, \em
Pure and Applied Mathematics, 115. Academic Press, Inc., Orlando, FL, 1984.}
\bibitem{[Ch]}\small{ S.S.\ Chern, \em On the curvatures of a piece of 
hypersurface in euclidean space. \em
 Abh.\ Math.\ Sem.\ Univ.\ Hamburg  {\bf 29} (1965), 77--91.}
\bibitem{[Ch2]}\small{ S.S.\ Chern, \em Simple proofs of two theorems 
on minimal surfaces. \em   Enseignement Math. (2)  15  1969 53--61.}
\bibitem{[D-H]}\small{ D.M.\ Duc; N.V.\ Hieu, \em
 Graphs with prescribed mean curvature on Poincar\'{e}  disk, \em
  Bull.\ London Math.\ Soc.\  {\bf 27} (1995),  no. 4, 353--358.}
\bibitem{[D-S]}\small{ D.M.\ Duc; I.M.C. Salavessa,
\em On a class of graphs with prescribed mean curvature, \em
 Manuscripta Math. {\bf 82}  (1994),  no. 3-4, 227--239.}
\bibitem{[F]}\small{ H.\ Flanders, \em  Remark on mean curvature, \em
 J.\ London Math.\ Soc.\  {\bf 41} (1966), 364--366.}
\bibitem{[Gaf]}\small{ M.P.\ Gaffney, \em A special Stokes's theorem 
for complete Riemannian manifolds.  Ann.\ of Math.\ (2)  {\bf 60}, 
(1954). 140--145.}
\bibitem{[Ger]}\small{ C.\ Gerhardt, \em $H$-surfaces in Lorentzian
 manifolds, \em  Comm.\ Math.\ Phys.\  {\bf 89} (1983),  no. 4, 523--553.}
\bibitem{[H]}\small{ E.\ Heinz, \em \"{U}ber Fl\"{a}chen mit eineindeutiger 
Projektion  auf eine Ebene, deren Kr\"{u}mmungen durch Ungleichungen 
eingeschr\"{a}nkt sind.\em 
Math.\ Ann.\  {\bf 129} (1955), 451--454.}
\bibitem{[N-R]}\small{B.\ Nelli; H.\ Rosenberg, \em 
Minimal surfaces in ${\mathbb H}\sp 2\times\mathbb R$, \em   
Bull.\ Braz.\ Math.\ 
Soc.\ (N.S.)  {\bf 33}  (2002),  no. 2, 263--292.}
\bibitem{[Ro]} \small{H.\ Rosenberg, \em  Minimal surfaces in 
${\mathbb M}\sp 2\times\mathbb R$. \em 
 Illinois J. Math. {\bf 46} (2002),  no. 4, 1177--1195.}
\bibitem{[S1]}
\small{I.M.C.\ Salavessa, \em
 Graphs with parallel mean curvature and a variational problem in
conformal geometry, \em Ph.D. Thesis, University of Warwick, 1987}
\bibitem{[S2]}
\small{I.M.C.\ Salavessa, \em Graphs with parallel mean curvature. \em
 Proc.\ Amer.\ Math.\ Soc.\  {\bf 107} (1989),  no. 2, 449--458.}
\end{thebibliography}
\end{document}